\documentclass{amsart}
\usepackage[utf8]{inputenc}

\title[Rationality of Nielsen zeta functions on solvmanfiolds]{On the rationality of the Nielsen zeta function for maps on solvmanifolds}
\author{Karel Dekimpe}\thanks{Research supported by long term structural funding - Methusalem grant of the Flemish Government.} 
\author{Iris Van den Bussche}

\newcommand{\End}[1]{\operatorname{End}(#1)}
\newcommand{\p}{\varphi}
\newcommand{\al}{\alpha}
\newcommand{\R}{\mathbb{R}}
\newcommand{\C}{\mathbb{C}}
\newcommand{\Q}{\mathbb{Q}}
\newcommand{\Z}{\mathbb{Z}}

\newcommand{\Aff}{\mathrm{Aff}}
\newcommand{\aff}{\mathrm{aff}}
\newcommand{\Aut}[1]{\mathrm{Aut}(#1)}

\newcommand{\GLZ}[1]{\mathrm{GL}_{#1}(\Z)}

\newcommand{\GLn}{\mathrm{GL_n(\Z)}}
\newcommand{\F}[1]{\operatorname{Fix}(#1)}

\newcommand{\ii}[1]{{#1}^{-1}}
\newcommand{\N}{\mathbb{N}}
\newcommand{\PI}{\Pi}
\newcommand{\tf}{torsion free }
\newcommand{\nor}{\triangleleft}
\newcommand{\NR}{$\mathcal{NR}$}

\newcommand{\infra}{infra-solv}
\renewcommand{\sf}{self-map }

\newcommand{\iso}[2]{\sqrt[#1\ ]{#2}}
\newcommand{\s}{strongly torsion free $S$-}
\newcommand{\di}[1]{\lvert\mathrm{det}(I-#1)\rvert}

\newcommand{\eig}[1]{\mathrm{eig}(#1)}
\newcommand{\til}{\tilde}

\newcommand{\vp}{polycyclic-by-finite }

\newcommand{\ad}{\mathrm{ad}}

\newcommand{\Lg}{{\mathfrak g}}
\newcommand{\sgn}[1]{\operatorname{sgn}(#1)}

\newcommand{\semn}{{\Z^n\rtimes_A\Z}}

\newcommand{\ZZ}{{\mathbb{Z}^2}}

\newcommand{\T}[1]{\mathrm{Tr}(#1)}

\newcommand{\simply}{connected, simply connected }

\newcommand{\stoms}{\Delta\backslash G}

\usepackage{tikz-cd}
\usepackage{amscd, amsfonts}
\usepackage[noadjust]{cite}
\usepackage{enumitem}
\usepackage{bm}
\usepackage{textcomp} % nice greek alphabet
\usepackage{pifont}   % Dingbats
\usepackage{booktabs}
\usepackage{amssymb,amsthm}
\usepackage{amsmath}
\usepackage{hyperref}

\newtheorem{theorem}{Theorem}[section]
\newtheorem{prop}[theorem]{Proposition}
\newtheorem{lemma}[theorem]{Lemma}

\newtheorem{con}[theorem]{Conjecture}
\theoremstyle{definition}
\newtheorem{definition}[theorem]{Definition}
\newtheorem{example}[theorem]{Example}
\newtheorem{claim}[theorem]{Claim}
\theoremstyle{remark}
\newtheorem{remark}[theorem]{Remark}

\begin{document}

\maketitle

\begin{abstract}
In \cite{dd15, fl15}, the Nielsen zeta function $N_f(z)$ has been shown to be rational if $f$ is a self-map of an infra-solvmanifold of type (R). It is, however, still unknown whether $N_f(z)$ is rational for self-maps on solvmanifolds.
In this paper, we prove that $N_f(z)$ is rational if $f$ is a self-map of a (compact) solvmanifold of dimension $\leq 5$. In any dimension, we show additionally that $N_f(z)$ is rational if $f$ is a self-map of an \NR-solvmanifold or a solvmanifold with fundamental group of the form $\Z^n\rtimes\Z$.
\end{abstract}

\section{Introduction}
Let $X$ be a compact manifold. 
Given a (continuous) self-map $f:X\to X$, we can define an integer called the Lefschetz number as
$$L(f)=\sum_{i=0}^{\mathrm{dim}(X)} (-1)^i\,\T{f_{\ast, i}: H_i(X,\Q)\to H_i(X,\Q)}.$$
The Lefschetz number is a homotopy invariant and indicates the presence of fixed points: if $L(f)\neq 0$, any map $g$ homotopic to $f$ has a fixed point. However, the Lefschetz number doesn't give us any information when $L(f)=0$, and it neither says something about the number of fixed points. 

Nielsen theory improves this setting as follows. Using the lifts of $f$ to the universal cover of $X$, one divides the fixed point set of $f$ into disjoint subsets called fixed point classes. Next, one assigns to each fixed point class an integer called the index of that fixed point class. If this index is nonzero, the fixed point class is called essential. 
The number of essential fixed point classes is then known as the Nielsen number of $f$; we denote it by $N(f)$. The Nielsen number is a homotopy invariant and in contrast with the Lefschetz number, $N(f)$ does give us some information about the number of fixed points: any map $g$ homotopic to $f$ has at least $N(f)$ fixed points. Moreover, Wecken showed in 1942 \cite{weck42-1} that if $X$ is not a surface, in fact
$$N(f)=\min\{ \#\F g \mid g \sim f\}\mbox{ \ where $\F g =\{ x \in X\mid g(x)=x\}$.}$$
Standard references for more information on the Lefschetz and Nielsen number include \cite{brow71-1} and \cite{jian83-1}.

In his 1967 paper \cite{smal67-1}, Smale introduced the Lefschetz zeta function of $f$ as the formal power series
$$L_f(z)=\exp\left(\sum_{k=1}^\infty \dfrac{L(f^k)}{k}z^k \right).$$
In that same paper, he showed that this function is in fact always a rational function.
Following Smale, Fel'sthyn introduced \cite{fels88-1, fp85-1} in 1985 the Nielsen zeta function of a self-map $f$ on a compact polyhedron $X$ as the formal power series 
$$N_f(z)=\exp\left(\sum_{k=1}^\infty \dfrac{N(f^k)}{k}z^k \right).$$
Unlike its Lefschetz counterpart, the Nielsen zeta function does not need to be rational in general.  
The rationality of the Nielsen zeta function has since been studied in several situations \cite{fels00-1, fels00-2, fh99-1, fp85-1, li94-1, wong01-1}.
Fairly recently, it was shown that $N_f(z)$ is rational if $f$ is a self-map of an infra-solvmanifold of type (R) \cite{dd15, fl15}. It is, however, still unknown whether $N_f(z)$ is rational for self-maps on all solvmanifolds.
In this paper, we settle this question in low dimensions:
\begin{theorem}
Let $f:S\to S$ be a \sf of a solvmanifold of dimension $\leq 5$. Then $N_f(z)$ is a rational function.
\end{theorem}
In the proof, we appeal numerous times to the related result mentioned above for Nielsen zeta functions on \infra manifolds of type (R) \cite{dd15, fl15}. %, Theorem~\ref{trar}.
Furthermore, we will use (and prove) the rationality of $N_f(z)$ in the following two cases:

\begin{prop}\label{nrar}
Let $f:S\to S$ be a \sf of an \NR-solvmanifold $S$. Then $N_f(z)$ is a rational function.
\end{prop}

\begin{prop}\label{znar}
Let $S$ be a solvmanifold with fundamental group $\PI=\Z^n\rtimes\Z$. Let $f:S\to S$ be a map. Then $N_f(z)$ is a rational function.
\end{prop}

After reviewing the necessary background in Section~\ref{prelim}, we prove Proposition~\ref{nrar} in Section~\ref{snr} and Proposition~\ref{znar} in Section~\ref{zn}. In Section~\ref{4}, we show that Propositions~\ref{nrar} and \ref{znar}
already cover all solvmanifolds up to dimension $4$. The remaining $5$-dimensional manifolds are then treated in Section~\ref{5}.
%Finally, in Section~\ref{6dim}, we formulate the conjecture that $N_f(z)$ is rational for all self-maps on solvmanifolds, and offer some partial results in dimension~6. 

%%%%%%%%%%%%%%%%%

\section{Solvmanifolds and infra-solvmanifolds of type (R)}\label{prelim}

Let $G$ be a \simply solvable Lie group. A \textit{solvmanifold} is a quotient space $\stoms$, with $\Delta$ a closed cocompact  subgroup. 
The fundamental group of a solvmanifold $\stoms$ is $K:=\Delta/\Delta_0$, with $\Delta_0$ the connected component of the identity in $\Delta$.
The algebraic structure of the fundamental groups of solvmanifolds is well known. Indeed,
by a result of Wang \cite{wang56-1}, a group $K$ is the fundamental group of a solvmanifold if and only if $K$ fits in a short exact sequence 
\begin{equation}\label{S-group}
1 \to N \to K \to \Z^k\to 1 
\end{equation}
where $N$ is a finitely generated torsion free nilpotent group. 
We will refer to a group $K$ fitting in a short exact sequence  of the above form (\ref{S-group}) as a \emph{\s group.}  
As $K/N\cong\Z^k$ is \tf abelian, we can always take $N=\iso{K}{[K,K]}$ in the above sequence, where for a group $G$ and a subgroup $H$ of $G$ we let $\iso{G}{H}$ denote the \textit{isolator} $\{g\in G\mid \exists k\in \N\setminus\{0\}: g^k\in H\}$ of $H$ in $G$. Hence, equivalently, a finitely generated group $K$ is a \s group if $\iso{K}{[K,K]}$ is finitely generated \tf nilpotent.

Let $M$ be a solvmanifold, and let $p:\tilde{M} \to M$ denote the universal covering projection, then  $\tilde{M}$ is diffeomorphic to $\R^n$ (with $n=\dim (M)$). The group of covering transformations of $p:\tilde{M} \to M$ is isomorphic to the fundamental group $K$ of $M$ and so we can view $M =  \tilde{M}/K$ (where $K$ is acting as the group of covering transformations).
Any self map $f:M \to M$ lifts to a map $\tilde{f}: \tilde{M} \to \tilde{M}$ (so with $p \circ \tilde{f} = f \circ p$). If we fix such a lifting $\tilde{f}$ (called the reference lift below), then there is an induced  endomorphism $\varphi: K \to K$ which is determined by 
\[ \forall k \in K, \forall \tilde{m} \in \tilde{M}:\; \tilde{f} ( k \cdot \tilde{m}) = \varphi(k) \cdot \tilde{f}(\tilde{m}).\]
Note that a different choice of reference lift will change $\varphi$ up to an inner autonorphism of $K$. Morever, up to a change of basepoint, the endomorphism $\varphi$ is exactly the endomorphism induced by $f$ on the fundamental group $K=\pi_1(M)$. 

\subsection{Model solvmanifolds and diagonal maps}

Model solvmanifolds and their corresponding diagonal maps were introduced by Heath and Keppelmann in \cite{hk02-1}. Their name is aptly chosen, as they form a class of relatively simple solvmanifolds which can be seen as `models' for Nielsen theory in the class of solvmanifolds. 
We collect here the relevant definitions and results from \cite{hk02-1}.

Let $K$ denote the semidirect product $\Z^n\rtimes_\al\Z^m$. Note that $K$ is a very simple \s group. Identifying $\R^{n+m}$ with $\R^n\oplus\R^m$, the group $K$ acts on $\R^{n+m}$ via
$$^{(z,k)}(x,y):=(z+\al(k)(x), k+y).$$

\begin{prop}\label{iseensolv}
Let  $K:=\Z^n\rtimes_\al\Z^m$. The resulting quotient space $\R^{n+m}/K$ is a solvmanifold with fundamental group $K$. 
\end{prop}
\begin{definition}
Let $K:=\Z^n\rtimes_\al\Z^m$, where the subgroup $\Z^n$ is fully invariant in $K$. The quotient space $\R^{n+m}/K$ from Proposition~\ref{iseensolv} is called a \textit{model solvmanifold}.
\end{definition}
%The minimal Mostow fibration of $M$ is given by $T^n\to M\to T^m$.

 We next define \textit{diagonal maps} on  model solvmanifolds. For $X\in M_n(\Z)$ and $Y\in M_m(\Z)$, let $(X,Y)$ denote the map $\R^{n+m}\to\R^{n+m}:(x,y)\mapsto (X(x), Y(y))$.
\begin{definition}
Let $M=\R^{n+m}/K$ be a model solvmanifold. A \sf $f:M\to M$ is said to be \textit{diagonal} if $f$ fits in a commutative diagram
$$
\begin{CD}
 \R^{n+m} @>(X,Y)>> & \R^{n+m}\\
 @ VVV & @ VVV\\
 M@>f>> & M
\end{CD}$$ for some $X\in M_n(\Z)$ and $Y\in M_m(\Z)$. We say that $f$ is diagonal of type $(X,Y)$ and for ease of notation, we also write $f=(X,Y)$.
\end{definition}
\begin{prop}
Let $M=\R^{n+m}/K$ be a model solvmanifold with fundamental group $K:=\Z^n\rtimes_\al\Z^m$.
\begin{enumerate}
    \item Let $f$ be a \sf of $M$ inducing the endomorphism $\p$ on $K$ and assume that $\p(\Z^n) \subseteq \Z^n$. Let  $X\in M_n(\Z)$ and $Y\in M_m(\Z)$ denote the induced endomorphisms on the subgroup $\Z^n$ and the factor $\Z^m$, respectively. Then $X\al(z)=\alpha(Yz)X$ for every $z\in\Z^m.$
    \item Given $X\in M_n(\Z)$ and $Y\in M_m(\Z)$ satisfying $X\al(z)=\alpha(Yz)X$ for every $z\in\Z^m,$ there exists a diagonal self-map $f$ of type $(X,Y)$. Moreover, $f$ induces the endomorphism $\p:K\to K:(z, k)\mapsto (Xz, Yk)$ with respect to the reference lift $(X,Y)$.
\end{enumerate}
\end{prop}

The following theorem says that for studying Nielsen theory on solvmanifolds, model solvmanifolds and diagonal maps are all we need.

\begin{theorem}
\label{aptlychosen}
Let $S$ be a solvmanifold. There exists a model solvmanifold $M$ (of the same dimension as $S$)  such that for every self-map $f:S\to S$, there exists a diagonal self-map $g:M\to M$ satisfying $N(f^k)=N(g^k)$ for every positive integer $k$. In particular, $N_f(z)=N_g(z)$.
\end{theorem}

\subsection{Infra-solvmanifolds of type (R)}\label{sectiontypeR}
Let $G$ be a \simply solvable Lie group. The affine group $\Aff(G)$ of $G$ is the semidirect product $\Aff(G)=G\rtimes \Aut{G}$. It embeds naturally in the semigroup $\aff(G)=G\rtimes\End{G}$. 
The product in $\aff(G)$ (and $\Aff(G)$) is given by $$(d,D)(e,E)= (dD(e), DE).$$
Both $\aff(G)$ and $\Aff(G)$ act on $G$ via $(d,D)\cdot g = d D(g)$.
One way to construct an \emph{infra-solvmanifold} is to consider a quotient manifold of the form $G/\Pi$ where 
$\Pi\subseteq \Aff(G)$ is a torsion free subgroup of the affine group of $G$ such that $\Gamma=G\cap \Pi$ is of finite index in $\Pi$ and 
$\Gamma$ is a discrete and cocompact subgroup of $G$.  
If the group $G$ is solvable of type (R) (also known as completely solvable), these types of infra-solvmanifolds $G/\Pi$ are called \emph{infra-solvmanifolds of type (R)}. Recall that a \simply solvable Lie group $G$ is said to be of type (R) if for every $X\in \Lg$ ($=\,$the corresponding Lie algebra of $G$)
the inner derivation $\ad(X)$ only has real eigenvalues. So for example, all nilpotent Lie groups are of type (R). 

Nielsen theory is well understood for infra-solvmanifolds of type (R). In particular, it can be shown \cite{dd15, fl15} that the Nielsen zeta function of every self-map of an infra-solvmanifold of type (R) is a rational function.

\section{Nielsen zeta functions on \texorpdfstring{\NR-solvmanifolds}{NR-solvmanifolds} are rational}\label{snr}

The class of \NR-solvmanifolds was introduced by Keppelmann and McCord in 1995 \cite{km95-1}  as a class of solvmanifolds satisfying the Anosov relation. %, that is, $N(f)=|L(f)|$ for every \sf $f$. 
We first recall the relevant properties of these manifolds in the following subsection.

\subsection{\texorpdfstring{\NR-solvmanifolds}{NR-solvmanifolds}.}
Let $K$ be a \s group and take $N=\iso{K}{[K,K]}$.
As $K$ is a \s group, $N$ is nilpotent, say of class $c$.
Let $\gamma_i(N)$ denote the $i$-th term of the lower central series of $N$, and put $N_i:=\iso{N}{\gamma_i(N)}$. 
Then $1\nor N_c\nor\cdots\nor N_1=N$ forms a central series of $N$ with free abelian factors $N_i/N_{i+1}$.
As the $N_i$ are normal subgroups of $K$, we get well-defined actions 
$$\rho_i:K/N\to \Aut{N_i/N_{i+1}}: \bar k\mapsto  \rho_i(\bar k) \mbox { with }\, \rho_i(\bar k):  xN_{i+1}\mapsto kx\ii k N_{i+1}$$ 
induced by conjugation. \label{rho_i}

\begin{definition}
We say that $K$ satisfies the \emph{\NR-property} if for every $i\in \{1,\ldots,c\}$ and for all $\bar k$ in $K/N$ the automorphism $\rho_i(\bar k)$ (of the free abelian group $N_i/N_{i+1}$) has no nontrivial roots of unity as eigenvalues.
\end{definition}
\begin{remark}
Here $\mathcal{NR}$ stands for ``No Roots''.
\end{remark}
\begin{example}\label{vbnr} The following are \NR-groups:
\begin{enumerate}
    \item If  $K=\Z^n$, then $N=1$, so plainly $\Z^n$ is $\mathcal{NR}$. %, as for every $z\in\Z^n/\{1\}$, the automorphism $\rho_1(z)=I_n$  has only eigenvalue $1$.
    \item For a less trivial example, consider the matrix 
$$A:=
\begin{pmatrix}
    2 & 1\\ 1 & 1
\end{pmatrix},$$
and let $K:=\Z^2\rtimes_\psi\Z$ with $\psi:\Z\to\GLZ2:t\mapsto A^t$.
Then $N=\Z^2$. 
It is easy to check that $A$ has only positive real eigenvalues. This of course holds equally for every power $A^t$. Hence $\rho_1(\bar t)=A^t$ has no nontrivial roots of unity as eigenvalues for every $\bar t\in K/\ZZ$, so $K$ is $\mathcal{NR}$.
\end{enumerate}
\end{example}

To see that the \NR-property does not really depend on the chosen series $1\nor N_c\nor \cdots\nor N_1=N$ of normal subgroups of $N$, we  introduce the following notations. 
Let $\p:\PI\to \PI$ be an endomorphism on a polycyclic-by-finite group $\PI$. 
Suppose that $$\PI_*: 1=\PI_{s+1}\nor \PI_s\nor \cdots \nor \PI_1=\PI$$ is a normal series of $\PI$ with finite or abelian factors $G_i:=\PI_i/\PI_{i+1}$ such that $\p(\PI_i)\subseteq \PI_i$ for every $i$ in $\{1,\ldots,s\}$. 
Then $\p$ induces endomorphisms $\p_i:G_i\to G_i$, which in turn induce endomorphisms $\bar \p_i$ on $G_i/\tau(G_i)$, where $\tau(G_i)$ is the set of torsion elements of $G_i$.
Note that $\tau(G_i)$ is indeed a subgroup of $G_i$ as $G_i$ is finite or abelian.
The groups $G_i/\tau(G_i)$ are free abelian groups of finite rank. 
Let $\eig {\bar\p_i}$ denote the set of eigenvalues of $\bar\p_i$, where we agree that $\eig{\bar\p_i}=\emptyset$ if $G_i/\tau(G_i)$ is trivial. 

\begin{lemma}
The set $\bigcup_{i=1}^s\eig{\bar\p_i}$ is independent of the chosen series.
\end{lemma}

This lemma can be proved by first showing that the set $\bigcup_{i=1}^s\eig{\bar\p_i}$ does not change if one refines the normal series and then by showing that two different normal series have ``equivalent'' refinements (See \cite[Theorem 8.4.3]{hall76-1}). 

\medskip
\label{alternativedefnr}
Accordingly, we will write $\eig \p:=\bigcup_{i=1}^s\eig{\bar\p_i}$. 
Using this notation, $K$ satisfies the $\mathcal{NR}$-property if and only if $\eig \p$ does not contain a nontrivial root of unity for every inner automorphism $\p$ of $K.$

\begin{definition}
A compact solvmanifold is an \emph{\NR-solvmanifold} if its fundamental group satisfies the \NR-property.
\end{definition}

Let $f:S\to S$ be a map on an \NR-solvmanifold $S$ with fundamental group $K$.  
Suppose that $f$ induces an endomorphism $f_*$ on $K$. Since $N=\iso{K}{[K,K]}$ is a fully characteristic subgroup of $K$, this endomorphism in turn induces an endomorphism $F_0$ on $K/N$ and endomorphisms $F_i$, $i=1,\ldots,c$, on the factor groups $N_i/N_{i+1}.$ 
The collection $\{F_0,\ldots, F_c\}$ is called the \emph{linearisation} of $f_*$.
Keppelmann and McCord  proved the following product formula for Nielsen numbers on \NR-solvmanifolds.

\begin{theorem}%[{\cite[Theorem~3.1]{km95-1}}]
[\cite{km95-1}]
\label{prod} 
Let $f:S\to S$ be a map on an \NR-solvmanifold $S$ with fundamental group $K$. 
Suppose that $f$ induces an endomorphism $f_*$ on $K$ with linearisation  $\{F_0,\ldots, F_c\}$. 
Then 
\begin{align*}
% &L(f)=\prod_{i=0}^c \,\det(I-F_i); \\%\quad\text{and}\quad
% &R(f)=\prod_{i=0}^c \,\di{F_i}_\infty;\\
&N(f)=\prod_{i=0}^c \,\di{F_i}.
\end{align*}
\end{theorem}

\subsection{The reduction technique}
We now prove that Nielsen zeta functions of self-maps on \NR-solvmanifolds are rational.
We will do so using what we call a \emph{reduction approach}:
\begin{itemize}
    \item we are given a self-map $f:S\to S$ of which we want to prove that $N_f(z)$ is rational;
    \item we construct an infra-solvmanifold $\tilde S$, of type (R), and a self-map $\tilde f:\tilde S\to\tilde S$ such that $N_{\tilde f}(z)=N_{f}(z)$; we say that $f$ \emph{reduces} to $\tilde f$;
    \item 
    we conclude that 
    %by Theorem~\ref{trar}, we conclude that 
    $N_{\tilde f}(z)=N_f(z)$ is rational. 
\end{itemize}

Using this approach, we can give a concise proof of Proposition~\ref{nrar}:
\begin{proof}[Proof of Proposition~\ref{nrar}]
Let $K$ be the fundamental group of $S$ and let $f_\ast$ be the endomorphism of $K$ induced by $f$ and consider the linearization $\{F_0, F_1, \ldots , F_c\}$ of $f_\ast$.
It is easy to see that $f^k$ induces $(f_\ast)^k$ which has linearization   $\{F_0^k, F_1^k, \ldots , F_c^k\}$ and so 
$N(f^k) = \prod_{i=0}^c\di{F_i^k}$.

Now consider $\Z^n$ with $n=k_0+k_1+ \cdots + k_n$ and let $F \in \mathrm{GL}_n(\Z)$ be the block diagonal matrix
$$F:=
\begin{pmatrix}
F_0 & \dots & 0\\
\vdots & \ddots& \vdots\\
0 & \dots & F_c
\end{pmatrix}.$$

Using this matrix $F$, we define the map \sf $\til f:\R^n/\Z^n\to\R^n/\Z^n: \bar z\mapsto \overline{Fz}$ on the torus which of course induces the endomorphism $F$ 
on $\Z^n$. Hence, $N(\til f)=\di{F}=N(f)$. 
Furthermore, $\widetilde{f^k}={(\til f\,)}^k$ for every $k\in\mathbb{N}$, hence, by the above reasoning, also $N(\til f^k)=N(f^k)$ for all $k\in\mathbb{N}$. So $N_f=N_{\til f}$. We have reduced $f$ to the self-map $\tilde f$.
As $N_{\til f}(z)$ is a rational function, we conclude that $N_f(z)$ is a rational function as well.
\end{proof}

%%%%%%%%%%%%%%

\section{Nielsen zeta functions on \texorpdfstring{$\Z^n\rtimes\Z$}{ZnxZ}}\label{zn}

Let $S$ be a model solvmanifold with fundamental group $\PI=\Z^n\rtimes\Z$. 
In this section, we prove Proposition~\ref{znar}: every self-map of $S$ has rational Nielsen zeta function. In the proof we again use the reduction approach, see Section~\ref{snr}. However, unlike in said section, we need some preliminary results before we can construct the self-map $\tilde f :\tilde S \to\tilde S$.

We first determine the Nielsen number of diagonal self-maps on $S$.
\begin{lemma}\label{form1}
Let $S$ be a model solvmanifold with fundamental group $\PI=\Z^n\rtimes_A\Z$. Let $f=(M,m)$ be a diagonal self-map of $S$. 
Then
$$N(f)=\frac{1}{d}\sum_{i=0}^{d-1} \lvert(1-m)\det(I-A^iM)\rvert$$
for any $d>0$ such that $A^d$ is $\mathcal{NR}$. 
\end{lemma}
\begin{proof}
Take $d>0$ such that $A^d$ is $\mathcal{NR}$ (it is clear that such a $d$ always exists). Then $K:=\Z^n\rtimes d\Z$ is $\mathcal{NR}$ and of finite index in $\PI$. 
Let $\p:\PI\to\PI:(z,t)\mapsto (Mz, mt)$ be the induced endomorphism by $f$. It is easy to see that $\p$ leaves $K$ invariant.
We can thus compute $N(f)$ using the averaging formula on infra-solvmanifolds, see \cite[Theorem~4.10 and Corollary~4.12]{dvnooit} %, Corollary~\ref{formc}.

%Set $K_1=\iso{K}{[K,K]}$.
Note that $K/\Z^n \cong\Z$ is \tf and abelian, so  $\iso{K}{[K,K]}\subseteq \Z^n$; say $\iso{K}{[K,K]} \cong \Z^s$ and $\Z^n\cong\Z^s\oplus\Z^r$. 
Then  $1 \nor \iso{K}{[K,K]} \nor K$ is a \tf filtration with free abelian factors 
$\iso{K}{[K,K]}\cong\Z^s$ and $K/\iso{K}{[K,K]}\cong\Z^r\times d\Z \cong \Z^{r+1}$.
We next determine the endomorphisms induced on these factors by the morphisms $\p$ and $\tau_i:\PI\to\PI: \gamma\mapsto (1,i)\gamma(1,-i)$.
%Thereto, we identify $\Z^n$ with $\Z^s\oplus\Z^r$ by choosing a basis according to this decomposition.

Choose a basis according to the decomposition $\Z^n=\Z^s\oplus\Z^r$. With respect to this basis, we can write
%As  $\Z^s=\iso{K}{[K,K]}$ is fully invariant, we can then write
 $$A=\begin{pmatrix}
    N& \ast \\
    0  & Q
    \end{pmatrix}\quad\text{and}
    \quad
    M=\begin{pmatrix}
    M_1& \ast \\
    0  & M_0
    \end{pmatrix}$$
for some $N\in\mathrm{GL}_s(\Z), Q\in\mathrm{GL}_r(\Z)$ and $M_1\in M_s(\Z), M_0\in M_r(\Z)$, %.
for $\Z^s=\iso{K}{[K,K]}$ is fully invariant.
Then $\tau_i$ and $\p$ induce the morphisms $N^i$ and $M_1$ on the first factor $\Z^s=\iso{K}{[K,K]}$, and the morphisms
$$
\begin{pmatrix}
         Q^i & 0 \\ 0 & 1
         \end{pmatrix}
\quad \text{and} \quad
\begin{pmatrix}
         M_0 & \ast \\ 0 & m
\end{pmatrix}
$$ on the second factor, respectively.
The averaging formula now asserts that
\begin{align*}
    N(f)
    &=\frac{1}{d}\sum_{i=0}^{d-1}\di{
    \left(\begin{smallmatrix}
         Q^i & 0 \\ 0 & 1
    \end{smallmatrix}\right)
    \left(\begin{smallmatrix}
         M_0 & \ast \\ 0 & m
\end{smallmatrix}\right)
    }
    \di{N^i M_1}\\
    &=\frac{1}{d}\sum_{i=0}^{d-1}\, \lvert(1-m)\det(I-Q^i M_0)\det(I-N^iM_1)\rvert\\
    &=\frac{1}{d}\sum_{i=0}^{d-1} \, \lvert (1-m)\det(I-A^i M)\rvert,
\end{align*}
completing the proof.
\end{proof}

The expression in Lemma~\ref{form1} closely resembles the averaging formula %\cite{ll09-1} 
for Nielsen numbers on \infra manifolds of type (R).
To exploit this resemblance, we decompose $A$ into a finite part and an $\mathcal{NR}$ part.
\begin{prop}\label{decom}
Let $A\in\GLn$, and let $d>0$ be minimal such that $A^d$ is $\mathcal{NR}$. Suppose $M\in M_n(\Z)$ and $m\in\Z\setminus\{0\}$ satisfy $MA=A^mM$.
Then $A$ has a decomposition $A=BC$ with $B, C\in\mathrm{GL}_n(\Q)$ satisfying
\begin{itemize}
    \item $C^d=I$ and $B$ is \NR
    \item $BC=CB$
    \item $MB=B^m M$ and $MC=C^m M$
\end{itemize}
\end{prop}

To prove this proposition, we record the following observation:
\begin{lemma}\label{sepeig}
    Let $S\in \mathrm{GL}_n(\Q)$ be diagonalisable. There exists $P\in \mathrm{GL}_n(\Q)$ such that
$$        
        PSP^{-1}=
        \begin{pmatrix}
                 S_b & 0\\ 0 & S_g
        \end{pmatrix}
$$
where $S_b\in\mathrm{GL}_k(\Q)$ has finite order and $S_g\in\mathrm{GL}_{n-k}(\Q)$ has no roots of unity as eigenvalues.   
    \end{lemma}
    
\begin{proof}
We appeal to the generalised Jordan canonical form \cite[Chapter~21, Theorem~5.4]{bjn94} of the matrix $S$: 
there exists $P\in\mathrm{GL}_n(\Q)$ and irreducible polynomials $p_1,\dots, p_s$ such that

 $$PSP^{-1}=
    \begin{pmatrix}
             \mathcal{J}_{r_1}(p_1) & &0\\
              &\ddots &\\
              0 & &\mathcal{J}_{r_s}(p_s)
    \end{pmatrix}.$$
Here $\mathcal{J}_{r}(p)$ is a generalised Jordan block
$$\mathcal{J}_r(p)=
      \begin{pmatrix}
                C(p) &\dots &\dots &0\\
                U &\ddots & &\vdots \\
                \vdots & \ddots &\ddots  &\vdots \\
                0 &\dots &U & C(p)
      \end{pmatrix}$$
built from the companion matrix of $p(t)=c_{0}+c_{1}t+\cdots +c_{{m-1}}t^{{m-1}}+t^m$
$$C(p):={\begin{pmatrix}0&0&\dots &0&-c_{0}\\1&0&\dots &0&-c_{1}\\0&1&\dots &0&-c_{2}\\\vdots &\vdots &\ddots &\vdots &\vdots \\0&0&\dots &1&-c_{{m-1}}\end{pmatrix}},$$
and the matrix
$$ U
      =\begin{pmatrix}
                0 &\dots &0 &1\\
                \vdots & &\vdots &0\\
                \vdots & &\vdots &\vdots\\
                0 &\dots &0 &0
      \end{pmatrix}.$$

Renumbering the polynomials $p_i$ if necessary, we may assume that $p_1,\dots,p_l$ are cyclotomic and $p_{l+1},\dots,p_s$ are not cyclotomic.
Set
$$        S_b:=
        \begin{pmatrix}
             \mathcal{J}_{r_1}(p_1) & &0\\
              &\ddots &\\
              0 & &\mathcal{J}_{r_l}(p_l)
        \end{pmatrix}\quad \text{and}\quad 
        S_g:=
        \begin{pmatrix}
             \mathcal{J}_{r_{l+1}}(p_{l+1}) & &0\\
              &\ddots &\\
              0 & &\mathcal{J}_{r_s}(p_s)
        \end{pmatrix}.$$
By construction $S_b$ has finite order, since $S$, hence $S_b$, is diagonalisable, and the eigenvalues of $S_b$, which are the roots of the cyclotomic polynomials $p_1,\dots, p_l$, are roots of unity.
In a similar vein, $S_g$ has no roots of unity as eigenvalues, since the $p_{l+i}$ are not cyclotomic.
\end{proof}

\begin{proof}[Proof of Proposition~\ref{decom}]
Let $A=US$ be the multiplicative Jordan decomposition of $A$. 
So $U$ is unipotent, $S$ is diagonalisable and $US=SU$.
As $A$ has rational entries, so has $S$ (and $U$). Hence, we can apply Lemma~\ref{sepeig} to find $P\in\mathrm{GL}_n(\Q)$ such that
$$         PSP^{-1}=
        \begin{pmatrix}
                 S_b & 0\\ 0 & S_g
        \end{pmatrix},
$$
where $S_b\in\mathrm{GL}_k(\Q)$ has finite order and $S_g\in\mathrm{GL}_{n-k}(\Q)$ has no roots of unity as eigenvalues. Note in particular that $S_b^{n_1}$ and $S_g^{n_2}$ have no eigenvalues in common for any $(n_1, n_2) \in\Z\times (\Z\setminus\{0\})$.
Set $$S_1:=P^{-1}
        \begin{pmatrix}
                 S_b & 0\\ 0 & I
        \end{pmatrix}P\quad\text{and}\quad S_2:=P^{-1}
        \begin{pmatrix}
                 I & 0\\ 0 & S_g
        \end{pmatrix}P.$$
We show:

\begin{claim}
Suppose $X\in M_n(\Q)$ and $x\in\Z\setminus\{0\}$ satisfy $XS=S^xX$. Then $XS_1=S_1^{x}X$ and $XS_2={S_2}^xX$.
\end{claim}
For a matrix $Y$, let $^PY$ denote $PYP^{-1}$. 
Put $^PM=\left(\begin{smallmatrix} \alpha & \beta \\ \gamma & \delta \end{smallmatrix}\right)$ with $\alpha\in M_k(\Q)$, $\beta\in\ \Q^{k\times (n-k)}$, $\gamma\in \Q^{(n-k)\times k}$ and $\delta\in M_{n-k}(\Q)$. Then
\vspace{-5pt}
\begin{align*}
MS=S^x M
&\Leftrightarrow {^PM}\ {^PS}= (^PS)^x\ ^PM\\
&\Leftrightarrow 
\begin{pmatrix} \alpha & \beta \\ \gamma & \delta \end{pmatrix}
\begin{pmatrix} S_b & 0 \\ 0 & S_g \end{pmatrix}
=
\begin{pmatrix} S_b^x & 0 \\ 0 & S_g^x \end{pmatrix}
\begin{pmatrix} \alpha & \beta \\ \gamma & \delta \end{pmatrix}\\
&\Leftrightarrow\begin{pmatrix} \alpha S_b\ & \beta S_g \\ \gamma S_b\ & \delta S_g \end{pmatrix}
=\begin{pmatrix} S_b^x \alpha\ & S_b^x \beta \\ S_g^x \gamma\ & S_g^x \delta \end{pmatrix}
\end{align*}
It is thus sufficient to prove that $\beta=\gamma=0$.
To this end, suppose first that $\beta\neq 0$. Then there exists $v\in\C^{n-k}$ with $S_gv=\lambda v$ and $\beta v\neq 0$. However, then $\lambda$ would be a common eigenvalue of $S_g$ and $S_b^x$, as $S_b^x\, \beta v=\beta\,S_g v=\beta\,\lambda v=\lambda\, \beta v$. Hence $\beta=0$.

Suppose next that $\gamma\neq 0$. Then there exists $v\in\C^{k}$ with $S_bv=\lambda v$ and $\gamma v\neq 0$. However, then $\lambda$ would be a common eigenvalue of $S_b$ and $S_g^x$, as $S_g^x\, \gamma v=\gamma\,S_b v=\gamma\,\lambda v=\lambda\, \gamma v$. Hence $\gamma=0$.
%Similarly $\gamma=0$.

\medskip

Consider $B:=US_2$ and $C:=S_1$. Note that $A=BC$ and $B, C\in\mathrm{GL}_n(\Q)$. We show that $B$ and $C$ satisfy the conditions of the proposition.
\begin{description}[font=\normalfont, itemsep=5pt]
    \item[$C^d=I$] By construction $S_b$ has only roots of unity as eigenvalues. As the eigenvalues of $S_b$ are also eigenvalues of $A$, and $A^d$ is \NR, all eigenvalues of $S_b$ are $d$-th roots of unity. As $S_b$ is diagonalisable, ${S_b}^d=I$, and thus also $C^d=I$.
    
    \item[$B$ is \NR] As $US=SU$, it follows at once from the claim above that $US_2=S_2U$. Thus $B$ is $\mathcal{NR}$ as $U$ is unipotent and $S_2$ is \NR.
    
    \item[$BC=CB$] Similarly, the claim above implies that $US_1=S_1U$. As $S_1S_2=S_2S_1$, too, $BC=CB$. 
    
    \item[$MB=B^mM$ and $MC=C^mM$] As $MA=A^mM$, also $MU=U^mM$ and $MS=S^mM$, see \cite[Lemma~4.3]{km95-1}. Again, the claim above implies that $MS_i={S_i}^mM$. Hence, $MC=C^mM$ and $MB=B^mM$  as well.
\end{description}
This completes the proof.
\end{proof}

Combining Lemma~\ref{form1} and Proposition~\ref{decom}, we can now prove Proposition~\ref{znar} using the same reduction approach as we used in Section~\ref{snr}.

\begin{proof}[Proof of Proposition~\ref{znar}]\label{hier}
Let $\PI:=\Z^n\rtimes_A\Z$ with $A\in\GLn$, and take $d>0$ minimal such that $A^d$ is \NR. Recall that we may assume that the self map $f$ is diganonal, so 
write $f=(M, m)$ with $M\in M_n(\Z)$ and $m\in\Z$. Then $MA=A^m M$.
By Lemma~\ref{form1}, $$N(f)=\frac{1}{d}\sum_{i=0}^{d-1} \lvert(1-m)\det(I-A^iM)\rvert.$$

We distinguish the cases $m=0$ (\textit{Case 1}) and $m\neq 0$ (\textit{Case 2}).

\noindent\textit{Case 1}. If $m=0$, then $MA=M$, hence, by induction, $MA^i=M$ for all positive integers $i$. Consider the map $\tilde f:\R^n/\Z^n\to\R^n/\Z^n:\bar x\mapsto \overline{Mx}$. Then
\begin{align*}
        N(f)&=\frac{1}{d}\sum_{i=0}^{d-1} \lvert\det(I-A^i M)\rvert\\
        &=\frac{1}{d}\sum_{i=0}^{d-1} \lvert \det(A^{-i})\det(I-A^i M)\det(A^i)\rvert\\
        &=\frac{1}{d}\sum_{i=0}^{d-1} \lvert \det(I- MA^i)\rvert\\
        &=\frac{1}{d}\sum_{i=0}^{d-1} \lvert \det(I- M)\rvert
        =N(\tilde f).
\end{align*}
Applying the above reasoning to $f^k$, $k\in\mathbb{N}$, also $N(f^k)=N(\tilde f^k)$. Hence $f$ reduces to the map $\tilde f$ (on a torus), so that $N_f(z)=N_{\tilde f}(z)$ is rational.

\noindent\textit{Case 2}. If $m\neq 0$, we  
apply Proposition~\ref{decom} to find a decomposition $A=BC$ with $B$ and $C$ satisfying the conditions in said proposition. 

It easily follows from $MB=B^mM$ that $MB^z=B^{mz}M$ for all $z\in\Z$.
As $B$ and $C$ commute, also $C^iMB^z=B^{mz}C^iM$ for every $i\in\mathbb{N}$. 
Hence, we can apply \cite[Theorem~4.1]{km95-1} on the matrix $C^iM\in M_n(\Q)$, the endomorphism $\Z\to\mathrm{GL}_n(\Q):z\mapsto B^z$ and the matrix $m\in\Z^{1\times 1}$ to find that $\det(I-C^iM)=\det(I-B^zC^iM)$ for all $z\in\Z$ if $m\neq 1.$ Note that Keppelmann and McCord prove \cite[Theorem~4.1]{km95-1} for integral matrices and an endomorphism $\Z^m\to\mathrm{SL}_n(\Z)$, but their proof carries over to the present situation.

%We can thus rewrite $N(f)$ as 
Combined, these facts show that
\allowdisplaybreaks\begin{align*}
        N(f)&=\frac{1}{d}\sum_{i=0}^{d-1} \lvert (1-m)\det(I-A^i M)\rvert\\
        &=\frac{1}{d}\sum_{i=0}^{d-1} \lvert (1-m)\det(I-B^iC^i M)\rvert\\
        &=\frac{1}{d}\sum_{i=0}^{d-1} \lvert (1-m)\det(I-C^i M)\rvert.
\end{align*}

Consider the subgroup $L$ of $\R^n$ defined by
    \begin{center}
          $L:=\Z^n+C\Z^n+\dots+C^{d-1}\Z^n$.
    \end{center}
As $C\in\mathrm{GL}_n(\Q)$, there exists $t\in\mathbb{N}$ with $\Z^n\leq L \leq \frac{1}{t}\Z^n$, thus $L$ is a lattice of $\R^n$.

Note that $C(L)\subseteq L$ as $C^d=I$, so we can define the group $\til \PI:=L\rtimes_C\Z$.
As $MC=C^m M$, also $M(L)\subseteq L$, so we can define the endomorphism $$\til \p:\til\PI\to\til\PI:(l,k) \mapsto (M(l),{mk})$$
of $\tilde \PI.$    

The group $\til\Pi=L\rtimes_c\Z$ embeds as 

$$\iota: \til\Pi\hookrightarrow \R^{n+1}\rtimes \mathrm{GL}_{n+1}(\R): (l,k) \mapsto \left[(l,k)\,,
\left(
\begin{smallmatrix}
C^k & 0 \\ 0 & 1
\end{smallmatrix}\right)\right].$$
Then
% \begin{itemize}%\setlength{\itemsep}{5pt}
%     \item 
    $\til \Pi\cap\R^{n+1}=L\oplus d\Z$ is a lattice of $\R^{n+1}$, and
    % \item 
    $\til\PI\cap\R^{n+1}$ has finite index in $\til\PI$,
% \end{itemize} 
so $\til\PI$ is a  Bieberbach group with holonomy group $\left\{\left(\begin{smallmatrix}
C^k & 0 \\ 0 & 1
\end{smallmatrix}\right)\mid k=0,\dots, d-1\right\}\cong\Z_d$.
On $\iota(\til\PI)$, the endomorphism $\til \p$ takes the form
$$\ \til \p:
\left[(l,k)\,,
\left(
\begin{smallmatrix}
C^k & 0 \\ 0 & 1
\end{smallmatrix}\right)\right]
\mapsto
\left[(M(l),mk)\,,
\left(
\begin{smallmatrix}
C^{mk} & 0 \\ 0 & 1
\end{smallmatrix}\right)\right].$$
% Thus $\tilde \p$ induces on $\tilde\PI\cap \R^{n+1}$ the endomorphism $\left(\begin{smallmatrix}
% M & 0 \\ 0 & m
% \end{smallmatrix}\right)$.
Setting $D=\left(\begin{smallmatrix}
M & 0 \\ 0 & m
\end{smallmatrix}\right)$, it is easily verified that for all $\gamma\in\iota(\tilde\PI)$,
$$(D\circ \gamma )(x)=(\tilde\phi(\gamma)\circ D)(x)\quad \text{ for every }x\in\R^{n+1},$$ where $\gamma=\left[(l,k)\,,
\left(
\begin{smallmatrix}
C^k & 0 \\ 0 & 1
\end{smallmatrix}\right)\right]$ maps $x=(r,s)\in\R^n\times\R$ %($r\in\R^n$, $s\in\Z$) 
to the element
$(l+C^k(r), k+s)$.
Hence $\tilde \p$ is induced by the affine map $\tilde f:=\overline{(0,D)}:\R^{n+1}/\tilde\PI\to\R^{n+1}/\tilde\PI:\bar x\mapsto \overline{Dx}$.
Here $\R^{n+1}/\tilde\PI$ is a flat manifold, which is a special case of an infra-solvmanifold of type (R), so we know that $N_{\til f}(z)$ is rational.

Using the averaging formula for Nielsen numbers \cite[Theorem~4.3]{ll09-1}, we compute that
\begin{align*}
    N(\til f) 
&=\frac{1}{d}\sum_{i=0}^{d-1}
\,\di{
\left(\begin{smallmatrix}
C^i & 0 \\ 0 & 1
\end{smallmatrix}\right)
\left(\begin{smallmatrix}
M & 0 \\ 0 & m
\end{smallmatrix}\right) }
\\
&=\frac{1}{d}\sum_{i=0}^{d-1}\,\lvert\det(I-C^iM)(1-m)\rvert\\
&=N(f).
\end{align*}

Moreover, $\widetilde{f^i}={\til f}^{\,i}$ for all $i\in\mathbb{N}$.
Thus $N({\til f}^{\,i})=N(\widetilde{ f^i})=N({f}^{i})$,
% $N_\p(z)
% =e^{\sum_{k=1}^\infty \dfrac{N(\p^k)}{k}z^k}
% =e^{\sum_{k=1}^\infty \dfrac{N({\til\p}^{\,k})}{k}z^k}
% =N_{\til\p}(z)$
so that $N_f(z)=N_{\til f}(z)$ is rational.
\end{proof}

\section{Solvmanifolds up to dimension 4}\label{4}

Let $f:M\to M$ be a self-map on a model solvmanifold $M$ of dimension $\leq 4$. In this section, we show that the Nielsen zeta function of $f$ is rational by showing that 
$M$ can be seen as an \NR-solvmanifold or an infra-solvmanifold of type (R). %, so the Nielsen zeta function of $f$ is rational.
%Hence, the Nielsen zeta function of $f$ -- and thus more generally of any self-map on a solvmanifold of dimension $\leq 4$ -- is rational.
Hence, more generally, 
%also
Nielsen zeta functions of self-maps on solvmanifolds of dimension $\leq 4$ are rational.
%, too.

Recall from page~\pageref{alternativedefnr} that equivalently, (a \s group) $K$ is \NR{} if $\eig{\tau_k}$ does not contain a nontrivial root of unity for every inner automorphism $\tau_k$ of $K$. Similarly,
\begin{definition}
We say that a \tf \vp group $\PI$ is \textit{of type (R) } if $\eig{\tau_\gamma}\subseteq \R^+$ for every inner automorphism $\tau_\gamma$ of $\PI$. We say that $\PI$ is virtually of type (R) if $\PI$ has a finite index subgroup of type (R).
\end{definition}

We make some elementary observations.
% \begin{definition}
We additionally say that $A\in\GLZ n$ is $\mathcal{NR}$ if $A$ does not have a nontrivial root of unity as eigenvalue, and that $A$ is \textit{of type (R)} if all the eigenvalues of $A$ are real and positive.
% \end{definition}

\begin{lemma}\label{anraknr}
If $A$ is \NR{} (resp. of type (R)), also $A^k$ is \NR{} (resp. of type (R)) for every integer $k$.
\end{lemma}
Hence,
\begin{lemma}\label{triviaal}
If $A\in\GLZ n$ is \NR{} (resp. of type (R)), also  $\semn$ is \NR{} (resp. of type (R)).
\end{lemma}

It is sufficient to show that the fundamental group of $M$ 
 is $\mathcal{NR}$ or virtually of type (R). 
Let 
$\Pi$ be the fundamental group of $M$. In order to check that $\PI$ is $\mathcal{NR}$ or virtually of type (R), we make three observations.

\begin{lemma}\label{lemma1}
    For every $A\in\mathrm{GL}_2(\Z)$, there exists an integer $k>0$ such that $A^k$ is of type (R).
\end{lemma}
\begin{proof}
If the eigenvalues of $A$ are real, we can take $k=2$. 
%Up to conjugacy, there are just six finite cyclic subgroups of $\GLZ2$ \cite[p.179]{newman}. Therefore,
If the eigenvalues of $A$ are not real, they are roots of unity of order $3$, $4$ or $6$ \cite[p.179]{newman}. So we can take $k=12$. 
\end{proof}
%Finally, lemma~\ref{lemma1} implies the following:
Lemma~\ref{lemma1} immediately implies the following:
\begin{lemma}\label{lemma3}
Let $m\in\{1,2\}$ and let $n$ be a positive integer. Then $\Z^m\rtimes \Z^n$ is virtually of type (R).
\end{lemma}

Lemma~\ref{lemma1} fails if $A\in\mathrm{GL}_3(\Z)$. We do have the following:

\begin{lemma}\label{lemma2}
     For every $A\in\mathrm{GL}_3(\Z)$ that is not $\mathcal{NR}$, there exists an integer $k>0$ such that $A^k$ is of type (R).
\end{lemma}

For future reference, we split the proof in two steps:
% This follows immediately from Lemma~\ref{nietnreigenwaardepm1} and: 
\begin{lemma}\label{nietnreigenwaardepm1}
     If $A\in\mathrm{GL}_3(\Z)$ is not $\mathcal{NR}$, then $A$ has eigenvalue $\pm1$.
\end{lemma}
\begin{proof}
If $A$ is not \NR, it has an eigenvalue $\lambda\neq 1$ which is a root of unity. If $\lambda=-1$, there is nothing to prove, so suppose $\lambda\neq -1$. Then $\lambda$ %\in\mathbb{C}\setminus \R$.
is not real, so $A$ has eigenvalues $\lambda$, $\overline\lambda$ and $\mu$ for some $\mu\in\R$. In fact, $\mu=\pm 1$, as 
$$1=\lvert\det(A)\rvert=|\lambda||\overline{\lambda}||\mu|=|\mu|$$
since $\lambda$ (thus also $\overline{\lambda}$) is a root of unity. 
\end{proof}

\begin{lemma}\label{eigenwaardepm1virttyper}
    If $A\in\mathrm{GL}_3(\Z)$ has $\pm1$ as an eigenvalue, there exists an integer $k>0$ such that $A^k$ is of type (R).
\end{lemma}
\begin{proof}
Let $\varepsilon\in\{-1,1\}$ be an eigenvalue of $A$. Take a corresponding eigenvector $v\in\Q^3$. Clearing denominators if necessary, we may assume that $v\in\Z^3$ and $\Z^3=\langle v , w, z\rangle$ for some $w, z\in\Z^3$. The matrix $A$ is thus similar over $\Z$ to a matrix of the form
$$\begin{pmatrix}
\varepsilon &\ast\\
0 & A'
\end{pmatrix}$$ with $A'\in\mathrm{GL}_2(\Z)$. Applying Lemma~\ref{lemma1}, we find $l>0$ such that ${A'}^l$ is of type (R), so $A^{2l}$ is of type (R) as well.
\end{proof}

We are now ready for

\begin{prop}
Let $S$ be a solvmanifold of dimension $\leq 4$. Then every self-map of $S$ has rational Nielsen zeta function.
\end{prop}
\begin{proof}
Let $f:S\to S$ be a self-map. Then there exists a model solvmanifold $\tilde S$, of dimension $\leq 4$, and a self-map $\tilde f$ of $\tilde S$ such that $N_f(z)=N_{\tilde f}(z).$ Let $\PI$ be the fundamental group of $\tilde S$. %Then $\PI=\Z^m\rtimes\Z^n$ with $m+n\leq 4$. 
We now list all the possibilities for $\PI=\Z^m\rtimes\Z^n$ and check that in each case, $\PI$ is $\mathcal{NR}$ or virtually of type (R).

\begin{itemize}
    \item In the $1$-, $2$- and $3$-dimensional case, $\PI\cong \Z$, $\PI\cong\Z^2$, $\PI\cong\Z\rtimes\Z$, $\PI\cong\Z^3$, $\PI\cong\Z\rtimes\Z^2$ or $\PI\cong\Z^2\rtimes\Z$ so $\PI$ is virtually of type (R) by Lemma~\ref{lemma3}.
    
    \item In the $4$-dimensional case, $\PI\cong\Z^4$, $\PI\cong\Z\rtimes\Z^3$, $\PI\cong\Z^2\rtimes\Z^2$ or $\PI\cong\Z^3\rtimes\Z$.
    %If $\PI\cong\Z\rtimes\Z^3$ or $\PI\cong\Z^2\rtimes\Z^2$, then
    In the former three cases,
    $\PI$ is again virtually of type (R) by Lemma~\ref{lemma3}.
    %If $\PI=\Z^3\rtimes_A\Z$
    So suppose that $\PI=\Z^3\rtimes_A\Z$. If $A$ is \NR, so is $\PI$. If $A$ is not \NR, Lemma~\ref{lemma2} delivers $k>0$ such that $A^k$ is of type (R). The finite index subgroup $\Z^3\times k\Z$ is then of type (R), hence $\PI$ is virtually of type (R) as well.
\end{itemize}
If $\PI$ is \NR, then $\tilde S$ is an \NR-solvmanifold, thus $N_f(z)=N_{\tilde f}(z)$ is rational by Proposition~\ref{nrar}.
If $\PI$ is virtually of type (R), then $\tilde S$ is an infra-solvmanifold of type (R). Hence $N_f(z)=N_{\tilde f}(z)$ is rational.  %Proposition~\ref{trar}.
\end{proof}

We conclude this section with an example of which Nielsen zeta functions can occur.

\begin{example} Let $S$ be a model solvmanifold with fundamental group $\PI=\Z^2\rtimes_A\Z$, where $A$ has real eigenvalues different from $\pm1$. We compute the possible Nielsen zeta functions on $S$.
So, let $f=(M,m)$ be a diagonal self-map of $S$.
Note that $S$ is an \NR-solvmanifold, so $$N(f^k)=\di{M^k}|1-m^k|$$ for every $k>0$.
We distinguish the cases %$m\neq\pm1$ (\textit{Case 1}), $m=1$ (\textit{Case 2}) and $m=-1$ (\textit{Case 3}).
$M=0$ (\textit{Case 1}) and $M\neq 0$ (\textit{Case 2}).

\noindent\textit{Case 1.}
If $M=0$, then $N(f^k)=|1-m^k|$.
We compute the Nielsen zeta function of $f$ similarly as in \cite{dtv18}. We thereto distinguish three cases:
\begin{itemize}
    \item If $m<-1$, then $|1-m^k|=(-1)^{k+1}(1-m^k)$. Hence 
    \begin{align*}
        N_f(z)&=\exp\left(\textstyle\sum_{k=1}^{+\infty} \, \frac{(-1)^{k+1}}{k}z^k - \frac{(-1)^{k+1}}{k}(mz)^k\right) \\
        &= e^{\ln(1+z) - \ln(1+mz)} \\
        &= \dfrac{1+z}{1+mz}.
    \end{align*}

    \item If $m\in\{-1,0,1\}$, then $|1-m^k|=1-m^k$, so similarly
    \begin{align*}
        N_f(z)&=\exp\left(\textstyle\sum_{k=1}^{+\infty} \, -\frac{(-1)^{k+1}}{k}(-z)^k + \frac{(-1)^{k+1}}{k}(-mz)^k\right) \\
        &= e^{-\ln(1-z) + \ln(1-mz)} \\
        &= \dfrac{1-mz}{1-z}.
    \end{align*}

    \item If $m>1$, then $|1-m^k|=-1+m^k$, hence
    \begin{align*}
        N_f(z)&= \dfrac{1-z}{1-mz}.
    \end{align*}
\end{itemize}
Moreover, for any $m$, the map $(0,m)$ defines a diagonal self-map on $S$, so all the above functions do arise as a Nielsen zeta function of a self-map on $S$.

\noindent\textit{Case 2.} 
Suppose next that $M\neq 0$. 
Then $m= \pm1$, for the condition $MA=A^mM$ only allows the trivial solution $M=0$ when $m\neq\pm1$. Indeed, suppose for a contradiction that $M\neq 0$. As $A$ is diagonalisable, there then exists $v\in\R^2$ such that $Mv\neq 0$ and $Av=\lambda v$. However, then also $A^m\, Mv=\lambda Mv$, hence $\lambda$ is an eigenvalue of $A^m$, too. This implies that $|\lambda|=|\lambda|^{\pm m}$, which is clearly impossible as $|\lambda|\neq1$ if $m\neq\pm1$. 
So $m=1$ or $m=-1$.

If $m=1$, clearly $N(f^k)=0$ for all $k>0$, hence $N_f(z)=1$. For instance, the identity self-map has Nielsen zeta function $1$.

\newcommand{\tr}[1]{\mathrm{Tr}(#1)}
So suppose that $m=-1$. Gon\c calves and Wong showed \cite{gw03-1} that %a nonzero 
$M$ in $M_2(\Z)$ satisfies $MA=A^{-1}M$ only if $\mathrm{Tr}(M)=0$. 
Moreover, for any $k>0$, the matrix $M^k$ satisfies $M^kA=A^{(-1)^k}M^k$, hence, if $k$ is odd,  $M^k$ satisfies $M^kA=A^{-1}M^k$ as well. 
Thus $\tr {M^k}=0$ whenever $k$ is odd. Note that if $k$ is even, $|1-m^k|=0$.
Hence, for any $k>0$,
\allowdisplaybreaks
\begin{align*}
    N(f^k)&=\di{M^k}|1-m^k|\\&=|1-m^k||1-\tr {M^k} +\det(M^k)|\\&=(1-(-1)^k))|1+\det(M)^k|.
\end{align*}
We again distinguish cases.
\begin{itemize}
    \item If $\det(M)<-1$, then $|1+\det(M)^k|=(-1)^k+(-\det(M))^k$, so $N(f^k)=(-1)^k+(-\det(M))^k-1-\det(M)^k$. Similarly as above, we can compute that
    \begin{align*}
        N_f(z)=\dfrac{(1-z)(1-\det(M)z)}{(1+z)(1+\det(M)z)} .
    \end{align*}
    \item If $\det(M)\geq-1$, then $1+\det(M)^k$ is positive, so $N(f^k)=1+\det(M)^k-(-1)^k-(-\det(M))^k$. Hence
    \begin{align*}
        N_f(z)=\dfrac{(1+z)(1+\det(M)z)} {(1-z)(1-\det(M)z)}.
    \end{align*}
     Note that if $\det(M)=0$, this agrees with what we found earlier in the case $M=0$ and $m=-1$.
\end{itemize}
This time, however, not every $\delta\in\Z$ will occur as the determinant of $M$ with $(M,-1)$ a diagonal self-map on $S$, or, equivalently, with $MA=A^{-1}M$.  Writing 
$$A=\begin{pmatrix}
 a & b\\ c & d
\end{pmatrix}
\quad\quad\text{and}\quad\quad
M=\begin{pmatrix}
 m & n\\p&-m
\end{pmatrix},$$
Gon\c calves and Wong showed in the same paper \cite{gw03-1} that $M$ satisfies $MA=A^{-1}M$ if and only if $(a-d)m+bp+cn=0$. As $A\neq\pm I_2$, this condition is not trivial, i.e., $(a-d,b,c)\neq(0,0,0)$. Hence, we can eliminate one variable ($m$, $n$ or $p$) and express $\det(M)$ as a homogeneous polynomial of degree 2 in the remaining two variables with rational coefficients. 

For example, take $A=\left(\begin{smallmatrix}
1 & 1\\ 1 & 2
\end{smallmatrix}\right)$.
Then $M=\left(\begin{smallmatrix}
m & n\\ p & -m
\end{smallmatrix}\right)$ satisfies $MA=A^{-1}M$ if and only if $m=n+p$. So all possible determinants of such $M$ are given in the set $\{-p^2-n^2-3np\mid n,p\in\Z\}$.
Hence, in this case, the possibilities for $N_f(z)$ are
\begin{itemize}
    \item $1$;
    \item $\dfrac{1+z}{1+mz}$ with $m\in\Z$ and $m<-1$;
    \item $\dfrac{1-z}{1-mz}$ with $m\in\Z$ and $m>1$;
    \item $\dfrac{1-mz}{1-z}$ with $m\in\Z$ and $|m|\leq 1$;
    \item ${\dfrac{(1+z)(1+\delta z)}{(1-z)(1-\delta z)}}^{\,\sgn {\delta+1}}$ with $\delta\in\{-p^2-n^2-3np\mid n,p\in\Z\}$.
\end{itemize}
\end{example}

\section{Solvmanifolds of dimension 5}\label{5}
Let $f:M\to M$ be a self-map on a model solvmanifold $M$ of dimension~$5$. In this section, we show that the Nielsen zeta function of $f$ is a rational function. Hence, the Nielsen zeta function of any self-map on a solvmanifold of dimension $5$ is rational as well.

Dimension $5$ is the smallest dimension in which (model) solvmanifolds exist that are neither \NR-solvmanifolds nor infra-solvmanifolds of type (R). 
Such manifolds, however, still are rare: in essence there are only two types of model solvmanifolds that meet this condition. We will first show that Nielsen zeta functions on these types of model solvmanifold are rational; next we show that all other model solvmanifolds are \NR-solvmanifolds or infra-solvmanifolds of type (R).

We start with a technical lemma:
\begin{lemma}\label{technischlemma}
Let $B\in\mathrm{GL}_3(\Z)$ and $M\in M_3(\Z)\setminus\{0\}$. Suppose that $B^k$ is not of type (R) for every $k>0$. 
\begin{itemize}
    \item If $MB=(-I)^\beta B^\delta M$ for some $\beta, \delta \in\Z$, then $\beta$ is even and $\delta=1$.
    \item If $B^\gamma M = (-I)^\xi M$ for some $\gamma, \xi \in\Z$, then $\xi$ is even and $\gamma=0$. 
\end{itemize}
\end{lemma}
\begin{proof}

As $B^k$ is not of type (R) for every $k>0$, the matrix $B$ has a complex eigenvalue $\omega$ such that $\omega^k$ is not real for every $k>0$. The eigenvalues of $B$ are then $\{\omega, \overline\omega,\lambda\}$ with $\lambda=\pm1/|\omega|^2\in\R$. 
Note that $B$ is \NR, see Lemma~\ref{lemma2}, hence $\lambda\neq -1$. Moreover,  by Lemma~\ref{eigenwaardepm1virttyper}, also $\lambda\neq1$ as $B^k$ is not of type (R) for every $k>0$. So every eigenvalue $\mu$ of $B$ has $|\mu|\neq 1$.

Suppose first that $B^\gamma M = (-I)^\xi M$.
As $M\neq 0$, there exists $v\in\Z^3$ with $Mv\neq 0$. Then
$$B^\gamma Mv= (-I)^{\xi}Mv = (-1)^{\xi}Mv,$$
so $(-1)^{\xi}$ is an eigenvalue of $B^\gamma$. As $B$ is $\mathcal{NR}$, this forces $\xi$ to be even, and as $B$ does not have $1$ as an eigenvalue, $\gamma$ must be zero.

Suppose next that $MB=(-I)^\beta B^\delta M$. 
The eigenvalues of $B$ are all distinct, so $B$ is diagonalisable. Since $M\neq 0$, there exists an eigenvector $v\in\mathbb{C}^3$ of $B$ such that $Mv\neq 0$. 
Let $\mu\in\{\omega, \overline\omega,\lambda\}$ be the eigenvalue corresponding to the eigenvector $v$. Using $MB=(-I)^\beta B^\delta M$, we see that
$$\mu Mv= MB v = (-I)^\beta B^\delta\, M v,$$
hence $\mu$ is an eigenvalue of $(-I)^\beta B^\delta$, too. There thus exists $\nu\in\{\omega, \overline\omega,\lambda\}$ satisfying $\mu=(-1)^\beta \nu^\delta$.
%Suppose now for a contradiction that $\beta$ is odd, so $\mu=- \nu^\delta$.
Switching $\omega$ and $\overline\omega$ if necessary, we may assume that $\nu=\lambda$ or $\nu=\omega$. This leaves the following possibilities. %for $\mu$ and $\nu$.
\begin{itemize}
    \item If $\nu=\lambda,$ also $\mu=\lambda$, then $\mu^{1-\delta}=(-1)^\beta$. As $\mu$ is not a root of unity, $\delta=1$. Hence $\beta$ is even. 
    
    \item If $\nu=\omega$ and $\mu=\omega$, again $\mu^{1-\delta}=(-1)^\beta$, implying $\delta=1$ and $\beta$ even. 
    
    \item If $\nu=\omega$ and $\mu=\overline\omega$, we find that  $\overline\omega=(-1)^\beta\omega^\delta$, hence $|\omega|=|\overline{\omega}|=|\omega|^\delta$. As $|\omega|\neq 1$, this implies $\delta=1$. Suppose for a contradiction that $\beta$ is odd. Then $\omega=-\overline\omega$, so $\omega=|\omega|i$ squared is real, a contradiction.
    
    \item If $\nu=\omega$ and $\mu=\lambda$, then $\omega^\delta=(-1)^\beta\lambda$,  contradicting that $\omega^\delta$ (or $\omega^{-\delta}$) is not real. So this case is impossible.
\end{itemize}
We conclude that $\beta$ must be even and $\delta=1$.
\end{proof}

\begin{prop}\label{speciaalgeval}
Let $S$ be a solvmanifold with fundamental group $\PI=\Z^3\rtimes_\psi\Z^2$. Suppose there are generators $t$ and $s$ of the factor $\Z^2$ such that 
\begin{itemize}
    \item $\psi(t)=-I$;
    \item $\psi(s)^k$ is not of type (R) for every $k>0$. 
\end{itemize}
Then every diagonal self-map of $S$ has a rational Nielsen zeta function.
\end{prop}
\begin{remark}
In the above proposition, $\psi(s)$ is necessarily \NR, see Lemma~\ref{lemma2}. 
\end{remark}
\begin{proof}
For ease of notation, let us set $B:=\psi(s)$.
Let $f:S\to S$ be a diagonal self-map of $S$. We compute the Nielsen number of $f$ using the averaging formula for Nielsen numbers on infra-solvmanifolds, \cite[Corollary~4.12]{dvnooit}. 

Write $f=(M,\Phi)$ for some $M\in M_3(\Z)$ and $\Phi\in M_2(\Z)$.
We distinguish the cases $M=0$ (\textit{Case 1}) and $M\neq 0$ (\textit{Case 2}). 

\medskip\noindent\textit{Case 1.} Suppose first that $M=0$. 
Consider the subgroup $K:=\Z^3\times \langle t^2, s^2\rangle$. Then $K$ is a normal, finite index, $\mathcal{NR}$ subgroup of $\PI$. Note that $f$ induces the endomorphism $\p:\PI\to\PI: (x,y)\mapsto (Mx, \Phi y)$, and that $\p$ leaves $K$ invariant.
Let $\tau_{t^is^j}$ denote conjugation by $t^is^j$ on $\PI$. 
The group $K$ has the torsionfree filtration $$1 \nor \Z^3 \nor K$$ with factors $\Z^3$ and $K/\Z^3=\langle t^2, s^2\rangle\cong \Z^2$. The induced endomorphisms by $\p$ on these factors are $M=0$ and $\Phi$, respectively; 
the endomorphism $\tau_{t^is^j}$ induces the endomorphism $(-I)^iB^j$ on the factor $\Z^3$ and the identity on the factor $K/\Z^3$. Hence, $f$ has Nielsen number
%\begin{align*}
$$    N(f)=\frac{1}{4}\sum_{i=0}^1\sum_{j=0}^1\, \lvert \det(I-(-I)^iB^jM)\det(I-\Phi)\rvert.$$
%     &=\frac{1}{4}\sum_{i=0}^1\sum_{j=0}^1\, \lvert\det(I-\Phi)\rvert
%     =\di{\Phi}.
% \end{align*}
However, as $M=0$, this simplifies to $N(f)=\di{\Phi}$. 
Consider the map $\tilde f:\R^2/\Z^2\to\R^2/\Z^2: \bar x\mapsto \overline{\Phi x}$. Then $N(f)=N(\tilde f)$.
Moreover, for any $k>0$, the map $f^k$ is a diagonal map $(0, \Phi^k)$, so also $N(f^k)=N(\tilde f^k)$ so that $N_f(z)=N_{\tilde f}(z)$ is rational. %So suppose henceforth that $M\neq 0$.

\medskip\noindent\textit{Case 2.} Suppose next that $M\neq 0$.
Consider the subgroup $K:=\Z^3\times \langle t^2, s\rangle$. Then $K$ is $\mathcal{NR}$ and of finite index in $\PI$. 
Moreover, we claim that $\p$ leaves $K$ invariant. Indeed, write $\Phi( t)= t^\alpha s ^\gamma$ and $\Phi( s)= t^\beta s ^\delta$ for some $\alpha, \beta, \gamma, \delta \in\Z$, so $\Phi\sim
\left(\begin{smallmatrix}
\alpha &\beta\\\gamma & \delta
\end{smallmatrix}\right)$. We have to show that $\beta$ is even.
We know that $M\psi(q)=\psi(\Phi( q))M$ for every $q\in\Z^2$. 
Taking $q=s$, we see in particular that $MB=(-I)^\beta B^\delta M$.
As $B^k$ is not of type (R) for every $k>0$, Lemma~\ref{technischlemma} says that $\beta$ is even and $\delta=1$.
Hence $\p(K)\subseteq K$.

Let $\tau_{t^i}$ denote conjugation by $t^i$ on $\PI$. 
The group $K$ has the torsionfree filtration $$1 \nor \Z^3 \nor K$$ with factors $\Z^3$ and $K/\Z^3\cong \Z^2$. The induced endomorphisms by $\p$ on these factors are $M$ and 
$
\left(\begin{smallmatrix}
\alpha & \beta/2 \\ 2\gamma & \delta
\end{smallmatrix}\right)
$,
respectively; 
the endomorphism $\tau_{t^i}$ induces the endomorphism $(-I)^i$ on the factor $\Z^3$ and the identity on the factor $K/\Z^3$. Hence, $f$ has Nielsen number
$$N(f)=\frac{1}{2}\sum_{i=0}^1\, \lvert \det(I-(-I)^iM)\det(I-\Phi)\rvert.$$

We now take $q=t$ in the relation $M\psi(q)=\psi(\Phi( q))M$ to find that $-M=(-I)^\alpha B^\gamma M$, or, alternatively, $B^\gamma M = (-I)^{\alpha+1}M$.
We again apply Lemma~\ref{technischlemma} to find that $\gamma=0$. So $\gamma=0$ and $\delta=1$, implying $\det(I-\Phi)=0$, hence $N(f)=0$.

Moreover, for any $k>0$, the self-map $f^k$ is diagonal of the form $(M^k, \Phi^k)$ with $\Phi^k\sim
\left(\begin{smallmatrix}
\ast &\ast\\0 & 1
\end{smallmatrix}\right)$ (even if $M^k$ happens to be zero). 
Hence $N(f^k)=0$ for all $k>0$, so $N_f(z)$ is rational.
\end{proof}

The second type of model solvmanifold that is not $\mathcal{NR}$ nor infra of type (R) has fundamental group $\Z^4\rtimes_A\Z$ where $A$ is not $\mathcal{NR}$ and $A$ has no power that is of type (R).
We thus get rationality on this second type of model solvmanifold from Proposition~
\ref{znar}. 

We next show that all other model solvmanifolds are \NR-solvmanifolds or infra-solvmanifolds of type (R).
%As in Section~\ref{4}, it is sufficient to show that the fundamental group is $\mathcal{NR}$ or virtually of type (R).
We already know from Lemma~\ref{lemma3} that the groups $\Z\rtimes\Z^4$ and $\Z^2\rtimes \Z^3$ are virtually of type (R). Moreover, if  $A$ is $\mathcal{NR}$ or $A^k$ is of type (R) for some $k>0$, the group $\Z^4\rtimes_A\Z$ is  $\mathcal{NR}$ or virtually  of type (R), too. So we are left to show that the groups $\Z^3\rtimes\Z^2$, other than those from Proposition~\ref{speciaalgeval}, are $\mathcal{NR}$ or virtually of type (R).

We first make some preliminary observations.
\begin{lemma}
 If $z\neq 0$ is an element of $\Z^m$, there exists a positive integer $d$ and a basis $\{x_1, \dots, x_m\}$ of $\Z^m$ such that $(x_1)^d=z$. 
\end{lemma}
\begin{lemma}\label{erbestaatgenerator}
If $\PI=\Z^n\rtimes_\psi\Z^m$ is not \NR, there exists $x_1\in\Z^m$ such that $\psi(x_1)$ is not \NR{} and $\{x_1\}$ extends to a set of generators $\{x_1, x_2, \dots, x_m\}$ of $\Z^m$.
\end{lemma}

\begin{lemma}\label{3bij2}
     Let $\PI=\Z^3\rtimes_\psi\Z^2$. Suppose 
     %$\PI$ does not satisfy the condition in Proposition~\ref{speciaalgeval}, that is,
     there are no generators $t$ and $s$ of the factor $\Z^2$ such that 
\begin{itemize}
    \item $\psi(t)=-I$;
    \item $\psi(s)^k$ is not of type (R) for every $k>0$.
\end{itemize}
Then $\Pi$ is $\mathcal{NR}$ or virtually of type (R).
\end{lemma}
\begin{proof}
Suppose $\PI$ is not \NR. Take generators $t$ and $s$ of the factor $\Z^2$. For ease of notation, let us set $A:=\psi(t)$ and $B:=\psi(s)$. By Lemma~\ref{erbestaatgenerator} above, we can assume that $A$ is not \NR. From Lemma~\ref{nietnreigenwaardepm1}, we infer that $A$ has an eigenvalue $\varepsilon\in\{\pm1\}$.

Consider the subgroup $W_\varepsilon:=\{z\in\Z^3\mid A(z)=\varepsilon z\}$ of $\Z^3$. Then $\Z^3$ decomposes as $\Z^3=W_\varepsilon\oplus\Z^d$ for some $d\in\{0,1,2\}$. Note that $B$ leaves $W_\varepsilon$ invariant. We can thus write 
$$
A=
\begin{pmatrix}
\varepsilon I & \ast\\ 0 & A'
\end{pmatrix}
\quad \text{and} \quad 
B=
\begin{pmatrix}
P & \ast \\ 0 & B'
\end{pmatrix},$$
where $A'$, $B'$ and $P$ are invertible, integral matrices.
We consider cases based on $d$.

\begin{itemize}
    \item[$d=2:$] In this case, $P=\pm1$ and $A', B' \in\mathrm{GL}_2(\Z)$. Lemma~\ref{lemma1} delivers $k>0$ and $l>0$ such that both $A'^k$ and $B'^l$ are of type (R), so $A^{2k}$ and $B^{2l}$ are of type (R), too. As $A$ and $B$ commute, this implies that the finite index subgroup $K:=\Z^3\times \langle t^{2k}, s^{2l}\rangle$ is of type (R), hence $\PI$ is virtually of type (R).
    
    \item[$d=1:$] In this case, $A', B'\in\{\pm1\}$ and $P\in\mathrm{GL}_2(\Z)$. Again, we get from Lemma~\ref{lemma1} $k>0$ such that $P^k$ is of type (R). Then $A^2$ and $B^{2k}$ are both of type (R), thus $\PI$ is virtually of type (R).
    
    \item[$d=0:$] In this case, $A=-I$. By assumption we then can find $k>0$ such that $B^k$ is of type (R). Since also $A^2=I$ is of type (R), $\PI$ is virtually of type (R).

    \end{itemize}
    We thus conclude that $\PI$ is $\mathcal{NR}$ or virtually of type (R).
\end{proof}
 
In the above proof, we showed the following statement, which we single out for future reference.
\begin{lemma}\label{nietidentiekecommuterendematrices}
     Let $A\neq \pm I_3$ be an element of $\mathrm{GL}_3(\Z)$ and suppose that $A$ has an eigenvalue in $\{\pm 1\}$. Then any $B\in \mathrm{GL}_3(\Z)$ commuting with $A$ has a power of type (R).
\end{lemma}

We summarise our findings in the following

\begin{prop}
Let $S$ be a solvmanifold of dimension $5$. Then every self-map of $S$ has rational Nielsen zeta function.
\end{prop}
\begin{proof}
Let $f:S\to S$ be a self-map. There exist a model solvmanifold $\tilde S$, of dimension $5$, and a diagonal self-map $\tilde f$ of $\tilde S$ such that $N_f(z)=N_{\tilde f}(z).$ Let $\PI$ be the fundamental group of $\tilde S$. Then $\PI=\Z^m\rtimes\Z^n$ with $m+n=5$. 
%We now list all the possibilities for $\PI=\Z^m\rtimes\Z^n$ and check that in each case, $\PI$ is $\mathcal{NR}$ or virtually of type (R).
This leaves the following possibilities:
\begin{itemize}
    \item If $\PI=\Z^5$, $\PI=\Z\rtimes\Z^4$ or $\PI=\Z^2\rtimes\Z^3$, then $\PI$ is virtually of type (R) by Lemma~\ref{lemma3}. Hence $N_{\tilde f}(z)$ is rational. 
    
    \item If $\PI=\Z^3\rtimes\Z^2$, then $N_{\tilde f}(z)$ is rational by Proposition~\ref{speciaalgeval} if $\PI$ satisfies the condition in said proposition; otherwise $N_{\tilde f}(z)$ is rational by Lemma~\ref{3bij2} and Proposition~\ref{nrar}. %and \ref{trar}.
    
    \item If $\PI=\Z^4\rtimes\Z$, then $N_{\tilde f}(z)$ is rational by Proposition~\ref{znar}.
\end{itemize}
We conclude that $N_f(z)=N_{\tilde f}(z)$ is rational as well.
\end{proof}

\section{What about higher dimensions?}
In \cite{phdIris} also all model solvmanifolds of dimension 6 were investigated using a case distinction on the isomorphism type of the fundamental group. For almost all cases we could show that the Nielsen zeta function is rational. The only types of model solvmanifolds in dimension 6 for which we were not yet able to prove the rationality are manifolds with a fundamental group of the form  $\Z^4\rtimes_\psi\Z^2$ such that
\begin{itemize}
    \item every $X\in\psi(\Z^2)$ is $\mathcal{NR}$ or has a power of type (R) and
    \item $\psi(\Z^2)=\langle A, B\rangle$ with $A^d=I$ for some $d>0$ and $\det(I-A^l)\neq 0$ if $0<l<d$.
\end{itemize}

The results of this paper, and those in dimension 6 from \cite{phdIris}  and the fact that the Nielsen zeta function is rational for several other families of solvmanifolds lead us to the following conjecture.

\begin{con}\label{conjecture}
Let $S$ be a solvmanifold. Then every self-map of $S$ has rational Nielsen zeta function.
\end{con} 

And we can even state the same about infra-solvmanifolds.

%\bibliographystyle{acm}
%\bibliography{allpapers}

\end{document}